\documentclass[11pt]{article}

\usepackage{amsmath, amsthm, amssymb, amscd}
\usepackage{graphicx}
\usepackage{bm}
\usepackage{tikz}
\usepackage[title]{appendix}

\usepackage{hyperref}
\usepackage{wrapfig}
\newcommand{\argmax}{\mathop{\rm arg\max}}

\numberwithin{equation}{section}
\newtheorem{thm}{Theorem}
\newtheorem{lem}{Lemma}
\newtheorem{rem}{Remark}
\newtheorem{prop}{Proposition}

\newcommand{\E}{\mathbb{E}}
\renewcommand{\P}{\mathbb{P}}

\newcommand{\sign}{\operatorname{sign}}

\newcommand{\mT}{\mathcal{T}}

\newcommand{\mP}{\mathcal{P}}

\def\limsup{\mathop{\overline{\rm lim}}}

\begin{document}

% "Title of the paper"
\title{On an estimator achieving the adaptive rate in nonparametric regression under $L^p$-loss for all $1\leq p\leq \infty$}

%\author{JSH}
\author{Johannes Schmidt-Hieber\footnote{University of Leiden, Niels Bohrweg 1, 2333 CA Leiden, Netherlands.\newline {\small {\em Email:} \texttt{schmidthieberaj@math.leidenuniv.nl}}}}

% \author{Johannes Schmidt-Hieber
% \vspace{0.1cm} \\
% {\em CREST-ENSAE,} \\ {\em 3, Avenue Pierre Larousse, 92240 Malakoff} \\ {\small {\em Email:} \texttt{Johannes.Schmidt.Hieber@ensae.fr}}}

% \footnotetext[1]{The research of Axel Munk and Johannes Schmidt-Hieber was supported by DFG Grant FOR 916 and GK 1023.}
% \footnotetext[2]{Author for correspondence, email:  \texttt{munk@math.uni-goettingen.de}}

%\date{\today}
\date{}
\maketitle

\begin{abstract}
Consider nonparametric function estimation under $L^p$-loss. The minimax rate for estimation of the regression function over a H\"older ball with smoothness index $\beta$ is $n^{-\beta/(2\beta+1)}$ if $1\leq p<\infty$ and $(n/\log n)^{-\beta/(2\beta+1)}$ if $p=\infty.$ There are many known procedures that either attain this rate for $p=\infty$ but are suboptimal by a $\log n$ factor in the case $p<\infty$ or the other way around.   In this article, we construct an estimator that simultaneously achieves the optimal rates under $L^p$-risk for all $1\leq p\leq \infty$ without prior knowledge of $\beta.$ In  contrast to classical wavelet thresholding methods that kill small empirical  wavelet coefficients and keep large ones, it is essential for simultaneous adaptation that on each resolution level, the largest empirical wavelet coefficients are truncated. This leads to a completely different point of view on wavelet thresholding. The crucial part in the construction of the estimator is the size of the truncation level which is linked to the unknown smoothness index. Although estimation of the smoothness index is known to be a difficult task, there is a data-driven choice of the truncation level that is sufficiently precise for our purpose.
\end{abstract}

% 62Gxx		Nonparametric inference
% %62G05   	Estimation
% %62G08   	Nonparametric Regression
% %62G20   	Asymptotic properties
% 

\newpage

\section{Introduction}
Suppose we observe $(Y_t)_{t\in [0,1]}$ with
\begin{align}
  dY_t= f(t) dt + n^{-1/2} dW_t, \quad t\in [0,1].
  \label{eq.GWNmodel}
\end{align}
Here, $f$ is the unknown regression function and $(W_t)_{t\geq 0}$ a Brownian motion. This is the Gaussian white noise model, which can be viewed as a continuous version of the classical nonparametric regression model, where we observe $f(\tfrac in)$ subject to some additive Gaussian noise (cf. Brown and Low \cite{brown1996}). Estimation of the function $f$ is a key problem in nonparametric statistics and has been studied extensively; for an overview see Tsybakov \cite{tsyb}. The estimation accuracy is typically evaluated under $L^p$-loss, with particular interest in the cases $p=2$ and $p=\infty.$ The $L^2$-norm coincides with the Kullback-Leibler divergence and is therefore the intrinsic loss function for the white noise model. The $L^\infty$-loss allows for uniform control over $\widehat f-f$ which is desirable for many applications and has a clear visual interpretation. In this work, we construct an estimator that achieves the minimax adaptive rate simultaneously for any $L^p$-loss with $p \in [1,\infty].$

It is well-known that term-by-term hard thresholding of the empirical wavelet coefficients leads to the adaptive rate $(n/\log n)^{-\beta/(2\beta+1)}$ for $\beta$-smooth signals under $L^\infty$-loss (for an overview on wavelet estimation, cf. Antoniadis \cite{ant}, Cai \cite{cai2012} or Johnstone \cite{john}). Hard thresholding is very stable but also conservative in the sense that only few coefficients are selected. Consequently, the reconstruction might miss smaller details of the true regression function, resulting in suboptimal convergence rates under $L^p$-loss with $p<\infty.$ This makes hard thresholding less appealing for applications. Finer details of the regression function can be recovered using block thresholding rules (cf. Hall et {\it al.} \cite{hal99}, Cai \cite{cai1999, cai}, Cavalier and Tsybakov \cite{cav01}). These methods take a block of wavelet coefficients and keep it if and only if $L^2$-signal is detected. Given a proper choice of the tuning parameters, block thresholding estimators achieve the 'clean' adaptive rate $n^{-\beta/(2\beta+1)}$ for $\beta$-smooth signals under $L^2$-risk. On the downside, the reconstructions show occasionally artificial spikes which are due to the large deviation behavior of the noise. As there are rather few such spikes, the block thresholding idea can be tuned to be simultaneously adaptive with respect to $L^2$-risk and squared pointwise risk (cf. Cai \cite{cai} and Cai and Silverman \cite{cai2}). Under $L^\infty$-loss, the spikes lead, however, to suboptimal $\log n$-factors as can be seen for instance from Theorem \ref{thm.neccesary_cond} below. An estimator that is simultaneously adaptive with respect to all $L^p$-losses $p\in [1,\infty],$ will inherit the good properties of both the term-by-term and the block thresholding in the sense that it finds more signal than the hard thresholding procedure and at the same time avoids large artificial spikes. 

Kernel smoothing with adaptive bandwidth choice and thresholding methods also achieve near optimal behavior if the parameter space is taken to be a general Besov ball. Nemirovski \cite{nem}, Section 3 constructs a method that achieves the adaptive rate up to a $\log n$-factor. Kerkyacharian and Picard \cite{ker00} show that hard thresholding leads to nearly adaptive methods over wide ranges of Besov spaces.  Similar results for block thresholding have been obtained by Chesneau \cite{che}.

One would like to consider simultaneous optimality from a more abstract point of view. Suppose that for a given estimation problem there are estimators that are  optimal with respect to two different loss functions. Does there exists then a procedure that is simultaneously optimal? If this is the case, is there a generic method to construct a simultaneously optimal estimator?  These questions are of practical importance in order to construct stable reconstruction  methods, but have not been addressed in this generality in the literature. In particular, a naive approach based on averaging estimators that are optimal with respect to one of the loss functions will not give simultaneously optimal procedures.

What makes simultaneous estimation so hard is the 'different bandwidth problem'. Suppose we have two loss functions corresponding to different minimax rates. Estimators achieving these rates will typically differ in their amount of smoothing/regularization  and therefore rely on different bandwidth/smoothness parameters. There is no straightforward way to merge these two estimators into one without being suboptimal for at least one of the two minimax rates. As an example consider the non-adaptive regression problem with known smoothness index $\beta$ and suppose we are interested in construction of an estimator that achieves the minimax $L^p$-risk for $p\in \{2,\infty\}.$ To obtain a rate optimal estimator under $L^2$-loss and $L^\infty$-loss, the bandwidth of a kernel estimator has to be chosen of order $h_{n,2}=n^{-1/(2\beta+1)}$ and $h_{n,\infty}=\big(\tfrac{\log n}{n}\big)^{1/(2\beta+1)},$ respectively. Using the bandwidth $h_{n,\infty}$ for $L^2$-loss or vice versa, $h_{n,2}$ for $L^\infty$-loss will result in additional $\log n$ factors in the rate. Since for a kernel estimator we have to fix one bandwidth $h_n,$ this suggests that, even if $\beta$ is known the class of kernel estimators is not flexible enough to allow for simultaneous minimax rates.

If the different bandwidth problem does not appear and both estimators rely on the same bandwidth parameter then simultaneous adaptation is typically possible and some results can be found in the literature. Consider for instance nonparametric estimation of the regression function $f$ under $L^2$-loss and (squared) pointwise loss. Then, for both problems $h_n=n^{-1/(2\beta+1)}$ is the optimal bandwidth and simultaneously adaptive estimators are known (cf. Cai \cite{cai}). A related problem is to estimate the regression function $f$ and its derivatives simultaneously. Rewriting this in terms of loss functions, this means that we want to achieve the optimal rate of convergence of $f$ under $L^2$-loss and Sobolev norm loss. In this case the nonparametric minimax rate for estimation of the $k$-th derivative under $L^2$-loss is $n^{-(\beta-k)/(2\beta+1)}$ and the optimal bandwidth is again $n^{-1/(2\beta+1)}$ for all $k.$ Efromovich \cite{efromovich1998} showed that even the exact asymptotic risk can be obtained simultaneously.

In order to achieve simultaneous adaptation under any $L^p$-loss, we have therefore to find a method that is unaffected by the different bandwidth problem. Rewritten in terms of wavelet coefficients, the different bandwidth problem says that on some resolution levels there is only negligible $L^\infty$-signal but possibly relevant $L^p$-signal for finite $p.$ To understand the problem, we first study simultaneous minimax estimation for fixed smoothness index. The main idea will be to truncate large empirical wavelet coefficients with a truncation level that explicitly depends on $\beta,$ large enough not to affect the $L^p$-risk but removing peaks in the reconstruction which otherwise would cause suboptimal $L^\infty$-loss convergence. Let us stress that the proposed wavelet truncation is the converse of classical thresholding. Instead of keeping large wavelet coefficients and killing small ones, we keep the small coefficients and truncate the largest coefficients on each resolution level. The method can easily be shown to be simultaneously minimax (cf. Theorem \ref{thm.non_adapt}). If  the smoothness index is unknown, we show that the truncation level cannot be estimated well enough to make the minimax result for fixed $\beta$ directly applicable in the adaptive case. Instead, we propose a truncation value that depends on the detected $L^p$-signal and prove that it has similar properties as in the fixed smoothness case. We also show that slight modifications of the estimator fail to be simultaneously adaptive implying that the proposed truncation of relatively large empirical wavelet coefficients is indispensable.

The paper is organized as follows. In Section \ref{sec.non_adapt}, we derive the wavelet truncation estimator that achieves the minimax rate of convergence with respect to any $L^p$-risk with $p\in [1,\infty]$ given knowledge of the H\"older smoothness $\beta.$ In Section \ref{sec.adapt}, a refined version of the truncation level is introduced which does not depend on $\beta.$ From that we can construct an estimator that is simultaneously adaptive. Proofs can be found in Section \ref{sec.proofs}.

{\it Notation:} The indicator function is denoted by $\mathbf{1}()$ and we write $\|\cdot \|_p$ for the $L^p$-norm on $[0,1].$

\section{Simultaneous estimation for known smoothness}
\label{sec.non_adapt}

Let $(\phi,\psi)$ be a scaling and wavelet function generating an orthogonal wavelet basis on $L^2[0,1].$ We assume that the wavelet basis is $s$-regular, that is, $\phi$ and $\psi$ have compact support and are $s$-times continuously differentiable. An explicit construction of such a wavelet bases on $[0,1]$ has been obtained by Cohen et {\it al.} \cite{coh}. Any $f\in L^2[0,1]$ can be expanded in $L^2$-sense as
\begin{align*}
  f =\sum_{k} c_{k}\phi_k + \sum_{j=0}^\infty\sum_{k\in I_j} d_{j,k} \psi_{j,k}
  =\sum_{j=-1}^\infty\sum_{k\in I_j} d_{j,k} \psi_{j,k} ,
\end{align*}  
with $\psi_{-1,k}:=\phi_k:=\phi(\cdot-k)$ and $d_{-1,k}:=c_k:=\int f(u)\phi_k(u)du$ as well as for $j\geq 0,$  $\psi_{j,k}:=2^{j/2}\psi(2^j\cdot-k)$ and $d_{j,k}:=\int f(u) \psi_{j,k}(u)du.$ The set $I_j$ consists of all indices $k$ on resolution level $j.$ In the following, we frequently make use of the fact that the cardinality of $I_j$ is of the order of $2^j,$ that is, $|I_j| \asymp  2^j.$ Given the Gaussian white noise model \eqref{eq.GWNmodel}, we can work with the empirical wavelet coefficients
\begin{align*}
  Y_{j,k}:= \int \psi_{j,k}(u) dY_u =d_{j,k} +\frac 1{\sqrt{n}} \epsilon_{j,k}, \quad \text{for} \ \  j\geq -1,k \in I_j,
\end{align*}
where $(\epsilon_{j,k})_{j,k}$ is an array of i.i.d. standard normal random variables.

As parameter space, we consider the wavelet representation of a H\"older ball with smoothness index $\beta,$
\begin{align*}
	\Theta(\beta,Q) := \big\{ f=\sum_{j=-1}^\infty\sum_{k\in I_j} d_{j,k} \psi_{j,k} \ : \ |d_{j,k}|\leq Q2^{-\frac j2(2\beta+1)}, \ \ \text{for all} \  j\geq -1,k\in I_j  \big\}.
\end{align*}
It is well-known that for $\beta<s,$ this space coincides with the classical H\"older ball, in the sense that there is a $Q'>0$ such that $\Theta(\beta,Q)$ contains all functions $f$ for which $f^{[\beta]}$ exists with $[\beta]:= \max\{u\in \mathbb{N} |u<\beta\}$ and $\sup_{x\in [0,1]} |f(x)|+ \sup_{x,y\in [0,1], x\neq y} \ |f^{[\beta]}(x)-f^{[\beta]}(y)|/|x-y|^{\beta-[\beta]} \leq Q'$ (cf. Cohen et {\it al.} \cite{coh}, Theorem 4.4).

Throughout the following, we make constantly use of the Besov space embedding
\begin{align}
	\big\|\sum_{j=-1}^\infty\sum_{k\in I_j} d_{j,k} \psi_{j,k} \big\|_p\leq C \sum_{j=-1}^\infty 2^{j(\frac 12 -\frac 1p)} \big\|(d_{j,k})_{k\in I_j}\big\|_{\ell^p}, \quad \text{for all} \ 1\leq p\leq \infty
	\label{eq.LP_embedding}
\end{align}
with $C$ a constant, independent of $p,$ and $\|\cdot\|_{\ell^p}$ the $\ell^p$-norm. The inequality can be directly verfied using the compact support of $\psi$ and the triangle inequality for the sum over $j.$

Choose integer sequences $(J_n)_n:=(J_n(\beta))_n$ and $(\widetilde J_n)_n:=(\widetilde J_n(\beta))_n$ such that $2^{J_n} \asymp n^{1/(2\beta+1)}$ and $2^{\widetilde J_n}\asymp (\tfrac n{\log n}\big)^{1/(2\beta+1)}.$ Inequality \eqref{eq.LP_embedding} together with the bound on the wavelet coefficients gives for finite $p,$
\begin{align*}
  \big\|\sum_{j=J_n+1}^\infty \sum_{k\in I_j} d_{j,k}\psi_{j,k}\big\|_p\lesssim n^{-\beta/(2\beta+1)} \quad \text{and} \ \ \big\|\sum_{j=\widetilde J_n+1}^\infty \sum_{k\in I_j} d_{j,k} \psi_{j,k}\big\|_\infty\lesssim \big(\tfrac{\log n}n\big)^{-\beta/(2\beta+1)}.
\end{align*}
Thus, as far as rates of convergence are concerned, there is no relevant signal on resolution levels above $J_n,$ for $p$ finite, and $\widetilde J_n,$ for $p=\infty.$ Consequently, the wavelet decomposition of $f$ has three different regimes $(I)-(III)$, which are for $p<\infty,$
\begin{align}
 f=\underbrace{\sum_{j=-1}^{\widetilde J_n} \sum_{k\in I_j} d_{j,k}\psi_{j,k}}_{(I): \ L^p\text{- and} \ L^\infty\text{-signal}}
  +\underbrace{\sum_{j=\widetilde J_n+1}^{J_n} \sum_{k\in I_j} d_{j,k}\psi_{j,k}}_{(II):\ L^p\text{- but no} \ L^\infty\text{-signal}}
  +\underbrace{\sum_{j=J_n+1}^{\infty} \sum_{k\in I_j} d_{j,k}\psi_{j,k}.}_{(III): \ \text{neither \ }L^p\text{- nor} \ L^\infty\text{-signal}}
  \label{eq.f_regimes}
\end{align}
For resolution levels $j\leq J_n,$ a reasonable choice is to take $Y_{j,k}$ as an estimator of $d_{j,k}.$ However, it might happen that the function $Y_{j,k}\psi_{j,k}$ itself is not in $\Theta(\beta,Q).$ Then, instead of estimating $d_{j,k}$ by $Y_{j,k},$ we can improve by projecting $Y_{j,k}$ on $[-Q2^{-\frac j2(2\beta+1)}, Q2^{-\frac j2(2\beta+1)}]$ (for the moment $\beta$ and $Q$ are assumed to be known). Consequently, the estimator for the scaling/wavelet coefficients $\widetilde d_{j,k}$ is given by
\begin{align}
  \widetilde d_{j,k}:= \sign(Y_{j,k}) \big(|Y_{j,k}|\wedge Q2^{-\frac j2(2\beta+1)}\big)
  \label{eq.djk_hat_adapt}
\end{align}
and the estimator for $f$ is
\begin{align}
  \widetilde f :=\sum_{j=-1}^{J_n}\sum_{k\in I_j} \widetilde d_{j,k} \psi_{j,k}.
  \label{eq.f_non_adapt_def}
\end{align}

\begin{thm}
\label{thm.non_adapt}
Work in model \eqref{eq.GWNmodel}. For any $\beta, Q\in (0,\infty)$ and any $p \in [1,\infty],$ the estimator $\widetilde f$ attains the minimax rate for $L^p$-loss over the parameter space $\Theta(\beta,Q).$
\end{thm}

This shows that wavelet truncation with truncation level $Q2^{-\frac j2(2\beta+1)}$ is enough to achieve simultaneous minimax rates. In the next section, we extend the result to unknown smoothness index.

\section{Simultaneous adaptation}
\label{sec.adapt}

By imitating the ideas from the previous section, we derive an estimator that achieves the adaptive rate for any $L^p$-risk. The crucial point in the construction is the projection on the interval $[-Q2^{-\frac j2(2\beta+1)},Q2^{-\frac j2(2\beta+1)}],$ for which knowledge of $\beta$ is required. 

One might be tempted to search for a data-driven procedure $t_j=t_j((Y_{j,k})_{j,k})$ which is independent of $\beta$ and estimates the boundary well enough in the sense that for a constant $Q'$ and with high probability, $Q2^{-\frac j2(2\beta+1)}\leq t_j\leq Q' 2^{-\frac j2(2\beta+1)}$ for all $j$ and all $\beta.$ Replacing the truncation level $Q2^{-\frac j2(2\beta+1)}$ by $t_j$ in \eqref{eq.djk_hat_adapt}, the resulting estimator would be simultaneously adaptive. Unfortunately, such a $t_j$ does not exist, as otherwise it would lead to a construction of honest, adaptive confidence bands. Indeed, we could take the intersection $[Y_{j,k}-2\sqrt{\log n/n}, Y_{j,k}+2\sqrt{\log n/n}] \cap [-t_j, t_j]$ as confidence set for the wavelet coefficient $d_{j,k}$ and then build from that a honest, adaptive confidence band for $f.$ Since honest, adaptive confidence bands do not exist, we have a contradiction and thus, there is no estimator $t_j$ with the imposed properties. A precise statement is given in the following theorem

\begin{thm}
\label{thm.non_exist}
Suppose that $\beta\in \{\beta_1, \beta_2\}$ with $\beta_1\neq \beta_2$ and $\beta_1,\beta_2>1/2.$ Let $J$ be such that $2^J \asymp n$ and consider a fixed $Q.$ There exists no data driven procedure $t_j=t_j((Y_{j,k})_{j,k})$ which is independent of $\beta$ and satisfies for some finite constant $Q',$
\begin{align*}
	Q2^{-\frac j2(2\beta_i+1)}\leq t_j\leq Q' 2^{-\frac j2(2\beta_i+1)}, \quad \text{for all} \ j\leq J, \ \text{for all} \ f\in \Theta(\beta_i, Q), \ i=1,2
\end{align*}
with probability tending to one.
\end{thm}

The proof can be found in the appendix. To derive the contradiction we only used the fact that honest, adaptive confidence bands do not exist. There is a stronger version of the non-existence result, which states that for two H\"older balls with different smoothness index $\beta_1< \beta_2,$ the best possible honest confidence band shrinks with the worst case rate $(n/\log n)^{-\beta_1/(2\beta_1+1)}$ on both spaces. This is due to the fact that there are sequences of functions for which we cannot test consistently to which of the two H\" older  spaces they belong to. Applied to estimation of the truncation level $t_j$ this shows that for some regression functions $f,$ we cannot even hope to be at least close to the target $Q2^{-\frac j2 (2\beta+1)}.$ 

What makes our situation different to the honest, adaptive confidence problem is the following. Consider again the simplified adaptation problem with only two different smoothness indices $\beta_1<\beta_2,$ say. Suppose that we have a signal for which we cannot test whether the regression function has smoothness $\beta_1$ or $\beta_2.$ For adaptive confidence bands, we have then to construct a band shrinking with the slower rate $(n/\log n)^{-\beta_1/(2\beta_1+1)},$ as we cannot exclude the possibility that the smoothness index is indeed $\beta_1.$ In contrast, for simultaneous adaption, we could try to work in such cases with the truncation level $t_j$ for the larger smoothness index $\beta_2.$ This induces an additional error if the true signal has smoothness $\beta_1.$ We will have to show then, that this error is of the correct order.

\subsection{Construction of estimator}

In order to construct an appropriate data-driven truncation level, we first need to introduce some notation. Define $\mP= \mP_n$ as the set of powers of two such that $p+1\leq \log n/(\log \log n)^2,$ that is,
\begin{align}
	\mP :=\big\{ p: p=2^r \ \text{with} \  r=0, 1,\ldots \ \text{and} \  p+1\leq \log n/(\log \log n)^2\ \big\}.
\end{align}
Set 
\begin{align}
	K_p := (6\sqrt{2} e +1)p^{p/2} \approx 24.07 p^{p/2}
    \label{eq.Kp_def}
\end{align}
and define for any positive integer $j,$ 
\begin{align*}
  D_{j,p}&:= \big\{S: S\subset I_j \ \text{and} \ \sum_{k\in S} |Y_{j,k}|^p \geq 2^pK_p n^{-p/2} |I_j| \big\}, \\
  L_{j,p}&:= \min_{S\in D_{j,p}} |S|,
\end{align*}
with $L_{j,p}:=\infty$ if $D_{j,p}$ is empty. The quantity $L_{j,p}$ can be viewed as the minimal number of (ordered) empirical wavelet coefficients, such that the inequality in $D_{j,p}$ holds. Furthermore, set 
\begin{align*}
	 t_{j,p} := \frac 2{\sqrt{n}} \Big(\frac{K_p|I_j|}{L_{j,p}-1}\Big)^{1/p}.
\end{align*}
Let $J_0$ be the smallest $j$ for which $|I_j|>n^{u_n},$ where $u_n=1/\log\log n.$ Notice that $n^{u_n}$ grows slower than any polynomial power of $n.$ It is therefore enough to truncate wavelet coefficients on resolution levels $j>J_0,$ only. Pick a  $J$ such that $2^J\asymp n.$ For $J_0< j\leq J,$ the estimator for the $(j,k)$-th wavelet coefficient is given by
\begin{align}
  \widehat d_{j,k} := \sign(Y_{j,k})(|Y_{j,k}| \wedge t_j), \ \ \text{with} \ \  	t_j :=
	\begin{cases}
	\infty &\text{if} \ \max_{k\in I_j} |Y_{j,k}| > 3\sqrt{\tfrac {\log n}n}, \\
	\max_{p\in \mP} t_{j,p}  &\text{otherwise}
	\end{cases}
  \label{eq.widehatdjk_def}
\end{align}
using the rule $1/0=+\infty$ and the {\it wavelet truncation estimator} is then
\begin{align}
  \widehat f :=\sum_{j=-1}^{J_0}\sum_{k\in I_j} Y_{j,k}\psi_{j,k}+\sum_{j=J_0+1}^{J}\sum_{k\in I_j} \widehat d_{j,k} \psi_{j,k}.
  \label{eq.f_adapt_def}
\end{align}

\begin{thm}
\label{thm.main}
Work in model \eqref{eq.GWNmodel}. For any $\beta, Q \in (0,\infty),$ there exists a constant $C,$ such that
  \begin{align*}
    \limsup_{n\rightarrow \infty}\ \sup_{1\leq p\leq \infty} \ \sup_{f\in \Theta(\beta,Q)} \big[\frac{n^{\beta/(2\beta+1)} }{\sqrt{p}} + \big(\frac{n}{\log n}\big)^{\beta/(2\beta+1)}\big]
    \E_f \big[ \|\widehat f-f\|_p\big]\leq C<\infty.
  \end{align*}
In particular, the estimator $\widehat f$ is simultaneously adaptive with respect to any $L^p$-loss with $1\leq p\leq \infty.$
\end{thm}

For fixed and finite $p,$ the $L^p$-risk is $\lesssim \sqrt{p} n^{-\beta/(2\beta+1)}.$ The fact that $E^{1/p}|\epsilon|^p \asymp \sqrt{p}$ for a standard normal random variable $\epsilon,$ explains the factor $\sqrt{p}$ in the rate, which can also be found for instance in  Nemirovskii \cite{nem}, Equation (1.49). The statement is uniform over $p,$ which allows to choose $p$ depending on $n.$ If $p\leq (\log n)^{2\beta/(2\beta+1)},$ then, we obtain the rate $\sqrt{p} n^{-\beta/(2\beta+1)}.$ For $p$ growing faster than $ (\log n)^{2\beta/(2\beta+1)},$ the rate becomes  $(n/\log n)^{-\beta/(2\beta+1)}$ which coincides with the $L^\infty$-risk.

The key argument in the proof of Theorem \ref{thm.main} are the following two properties of the truncation level. 

\begin{prop}
\label{prop.properties_of_adapt_est}
{\ }{}
\begin{itemize}
\item[(i)] If $\max_{k\in I_j} |d_{j,k}|\leq \sqrt{\tfrac{\log n}n}$ and $f\in \Theta(\beta,Q),$ then, with probability $1-O(\log^2 n/n),$ 
\begin{align*}
	t_j\leq 4 Q 2^{-\frac j 2 (2\beta+1)}, \quad \text{for all} \ J_0 < j \leq J.
\end{align*}
\item[(ii)] $|\{ k: |Y_{j,k}| > t_j\}| < \min_{p\in \mP} L_{j,p}.$ 
\end{itemize}
\end{prop}

Recall that Theorem \ref{thm.non_exist} states that it is impossible to have a purely data-driven $t_j,$ such that $Q2^{-\frac j 2 (2\beta+1)}\leq t_j\leq Q' 2^{-\frac j 2 (2\beta+1)}$ with probability tending to one. The first statement of the proposition says that the truncation level $t_j$ satisfies the upper bound with $Q'=4Q.$ The second assertion of the proposition states that the number of truncated wavelet coefficients on resolution level $j$ is always bounded by the minimum of all $L_{j,p},$ $p\in \mP.$ This property implicitly controls the lower bound of $t_j$ and assures that the $L^p$-risk of the coefficient vector $(d_{j,k}\mathbf{1}(|Y_{j,k}| > t_j))_{j,k}$ is of the order $\sqrt{p} n^{-\beta/(2\beta+1)}$ (cf. part (IV) in the proof of Theorem \ref{thm.main}). Thus, truncation by $t_j$ does not affect the convergence rate under $L^p$-risk.

In \eqref{eq.widehatdjk_def}, we set $\widehat d_{j,k}=\sign(Y_{j,k}) t_j$ if $|Y_{j,k}|>t_j.$ Theorem \ref{thm.main} remains true for any choice with $|\widehat d_{j,k}|\leq t_j.$ Nevertheless, the truncation step \eqref{eq.widehatdjk_def} is necessary. To see this, consider the modified estimator which does not truncate empirical wavelet coefficients exceeding in absolute value the truncation level $t_j$ on resolution levels where $L^p$-signal was detected, that is, if $t_j>0.$ This estimator can be written as
\begin{align}
  \widehat f_1 := \sum_{j=-1}^{J_0}\sum_{k\in I_j} Y_{j,k} \psi_{j,k} + \sum_{j=J_0+1}^{J} \mathbf{1}(t_j>0) \sum_{k\in I_j} Y_{j,k} \psi_{j,k}.
  \label{eq.f_block_thresh}
\end{align}
However, as stated in the following theorem, $\widehat f_1$ is suboptimal by a $(\log n)^{1/(4\beta+2)}$-factor implying that truncation is indeed necessary.

\begin{thm}
\label{thm.neccesary_cond}
Work in model \eqref{eq.GWNmodel}. For any $\beta, Q \in (0,\infty),$
\begin{align*}
  \sup_{f\in \Theta(\beta,Q)} \E_f\big[\|\widehat f_1-f\|_{\infty}\big]
  \gtrsim \big(\frac{n}{\log n}\big)^{-\beta/(2\beta+1)} (\log n)^{1/(4\beta+2)}.
\end{align*}
Consequently, $\widehat f_1$ does not achieve the adaptive rate with respect to $L^\infty$-loss.
\end{thm}

The proof of Theorem \ref{thm.neccesary_cond} is based on the following observation: Consider $f\in \Theta(\beta,Q)$ with wavelet coefficients $d_{j,k}= Q2^{-\frac j2(2\beta+1)}$ for all $j,k.$ Then, for the thresholding estimator \eqref{eq.f_block_thresh}, we can find an integer $J_n=J_n(\beta)$ with $2^{J_n}\asymp n^{1/(2\beta+1)}$ such that $\widehat d_{j,k}=Y_{j,k}$ for $j\leq J_n$ and all $k.$ By the extreme value behavior of the maximum of standard normal random variables there exists, with high probability, a $k^*$ such that $\epsilon_{J_n,k^*}\geq q\sqrt{\log n}$ for some $q>0.$ Such a $k^*$ will lead to an artificial spike in the sense that
\begin{align}
  \|Y_{J_n,k^*}\psi_{J_n,k^*}\|_\infty
  \gtrsim \big(\frac{n}{\log n}\big)^{-\beta/(2\beta+1)}(\log n)^{1/(4\beta+2)}.
  \label{eq.lowerbd_heur}
\end{align}
By some refined analysis, we obtain that with high probability, this term is a lower bound of the $L^\infty$-loss, that is, $\|\widehat f_1-f\|_\infty\gtrsim \|Y_{J_n,k^*}\psi_{J_n,k^*}\|_\infty.$ Taking expectation and using \eqref{eq.lowerbd_heur}, the lower bound of Theorem \ref{thm.neccesary_cond} follows. A complete proof can be found in the appendix.

Simultaneous adaptation can also not be achieved by blockwise soft thresholding. In contrast to hard thresholding, a soft thresholding procedure shrinks large empirical wavelet coefficients as well and has thus a similar effect than wavelet truncation. It shrinks, however, also small coefficients by the same factor leading to suboptimal behavior. Let us describe this in more detail for known smoothness $\beta.$ Recall the decomposition \eqref{eq.f_regimes}. For any critical resolution level $j\in \{\widetilde J_n+1,\ldots,J_n\},$ let $a_j$ be a shrinkage factor and $\widehat{d}_{j,k} = a_{j,k}Y_{j,k}$ be the soft thresholding estimator for the wavelet coefficient $d_{j,k}.$ The maximum $\max_k |\epsilon_{j,k}|$ is of the order $\sqrt{\log n}$ and thus, in order to be adaptive under $L^\infty$-loss we must have $a_j\lesssim 2^{-\frac j2 (2\beta+1)}/\sqrt{\log n/n}.$ There exists $j^* \in  \{\widetilde J_n+1,\ldots,J_n\},$ such that for $n\rightarrow\infty,$ $2^{\widetilde J_n}\ll 2^{j^*} \ll 2^{J_n}.$ Consequently, $a_{j^*}\rightarrow 0.$ Consider a regression function $f$ for which $|d_{j,k}| \asymp 2^{-\frac j2 (2\beta+1)}.$  The $L^p$-signal on the $j^*$-th resolution level is $2^{j^*(\frac 12 -\frac 1p)} (\sum_k |d_{j^*,k}|^p)^{1/p} \gg n^{-\beta/(2\beta+1)}.$ Since $a_{j^*}\rightarrow 0,$ the $L^p$-bias of the reconstruction on the $j^*$-th resolution level will thus be of larger order than $n^{-\beta/(2\beta+1)},$ leading to a suboptimal rate of the soft thresholding estimator under $L^p$-loss. Hence, for simultaneous adaptation, wavelet truncation is necessary and classical thresholding methods do not work.

From a theoretical point of view, it remains unclear, whether the results could be extended to more general parameter spaces. Our proofs rely heavily on the assumption that the parameter space is a H\"older ball. Otherwise problems occur on high resolution levels, as we have to search for all $L^p$-signals with $p\leq p_n$ and $p_n$ slowly growing to infinity. To explain this in more detail, suppose that the parameter space is a $L^2$-Sobolev ball and that we are interested in the rate under $L^2$-loss. Then, it might well happen that even on large resolution levels $j,$ $L^p$-signal is detected with $p>2,$ which would then give $t_j\geq t_{j,p}>0.$ Thus, far too many empirical wavelet coefficients are included into the reconstruction leading to suboptimal $L^2$-risk. Notice that this problem only occurs if we want to adapt to all $L^p$-losses. If we were instead interested in simultaneous adaptation with respect to two loss functions, say $L^2$- and $L^\infty$-loss, a wavelet truncation estimator could be constructed that is simultaneously adaptive over Besov balls using a blockwise truncation scheme.

Let us shortly comment on applicability of the proposed estimator for real data. The aim of this work is to present the phenomena underlying simultaneous adaptation and in particular to introduce wavelet truncation as a tool to robustify estimators avoiding as many technicalities as possible. One can easily extend the results to larger classes of estimators that also achieve in theory the adaptive rate under any $L^p$-loss and lead to procedures with better finite sample properties. Instead of defining one global truncation level on each resolution level, for instance, there should be some gain in splitting the coefficients into smaller blocks and working locally. The performance of the refined estimator will then rely in practice on a careful choice of the block length. To improve on the applicability, one might want moreover to optimize the constant in the exponential inequality stated in  Lemma \ref{lem.exp_ineq}. This would then allow for a better constant in the definition of $D_{j,p}.$ To conclude the discussion on applicability, let us mention that wavelet truncation is computationally feasible. Indeed, sorting the empirical wavelet coefficients $Y_{j,k}$ on each resolution level according to their absolute value in a first step allows to calculate $t_j$ without explicitly computing $D_{j,p}.$

{\bf Acknowledgments.}
This research was partially supported by DFG postdoctoral fellowship SCHM 2807/1-1. The author is grateful for helpful discussions with Marc Hoffmann and Judith Rousseau during the preparation of a first version. He would further like to thank the participants of the Oberwolfach workshop 1411 for many interesting remarks and suggestions.

\section{Appendix}
\label{sec.proofs}

\begin{proof}[Proof of Theorem \ref{thm.non_exist}]
Suppose such a $t_j$ exists. Define the set $A= \{f =\sum_{j=-1}^\infty\sum_{k\in I_j} d_{j,k} \psi_{j,k} : d_{j,k} \in [Y_{j,k}-\sqrt{2\log n/n}, Y_{j,k}+\sqrt{2\log n/n}] \cap [-t_j, t_j], \ \text{if} \ j\leq J \ \text{and} \  |d_{j,k}|\leq  2^{-j} \log n, \ \text{if} \ j>J\}.$ We show that $A$ is a honest, adaptive confidence band in the sense that for $i=1,2,$ $\inf_{f\in \Theta(\beta_i, Q)} \P_f\big( f\in A) \rightarrow 1$ and $\sup_{f\in \Theta(\beta_i,Q)}\E_f \sup_{g, h\in A}\|g-h\|_\infty\lesssim n^{-\beta_i/(2\beta_i+1)}.$

Assume that $f\in \Theta(\beta, Q)$ for $\beta \in \{\beta_1,\beta_2\}.$ By construction $|d_{j,k}|\leq t_j$ for $j\leq J$ and for sufficiently large $n$ also $|d_{j,k}|\leq Q2^{-\frac j2 (2\beta+1)}\leq 2^{-j}\log n$ thanks to $\beta>1/2.$ Since the probability of $\max_{j\leq J,k\in I_j} |\epsilon_{j,k}|>\sqrt{2\log n}$ converges to zero as $n\rightarrow \infty,$ we obtain that $\inf_{f\in \Theta(\beta, Q)} \P_f\big( f\in A) \rightarrow 1.$  

To bound the diameter of the confidence band, recall $2^{\widetilde J_n(\beta)} \asymp (n/\log n)^{1/(2\beta+1)}$ and observe that with \eqref{eq.LP_embedding}, for any $g, h \in A,$
\begin{align*}
	\|g-h\|_{\infty} 
	&\leq 2C \sum_{j=-1}^{\widetilde J_n(\beta)} 2^{j/2}  \sqrt{2\log n/n} + 2C\sum_{j= \widetilde J_n(\beta)+1}^J 2^{j/2}  Q' 2^{-\frac j2 (2\beta+1)}  + 2C\sum_{j= J+1}^\infty 2^{j/2} 2^{-j} \log n \\
	&\lesssim (n/\log n)^{-\beta/(2\beta+1)}.
\end{align*}
Thus, $A$ is a honest, adaptive confidence band for the regression function. However, they do not exist as shown by Low \cite{low1997} for density estimation. The non-existence carries over to nonparametric regression, cf. Genovese and Wasserman \cite{genovese2008}. Therefore, we have a contradiction and such a sequence of $t_j$ cannot exist. This completes the proof.
\end{proof}

\begin{lem}
\label{lem.mom_norm}
Let $\epsilon\sim \mathcal{N}(0,1).$ For any $q\geq 1,$ $\E|\epsilon|^q \leq q^{q/2}.$
\end{lem}

\begin{proof}
For $r>2,$ partial integration gives $\E[|\epsilon|^r]=(r-1)\E[|\epsilon|^{r-2}].$ For $0<r\leq 2,$ Jensen's inequality shows that $\E[|\epsilon|^r]\leq  \E^{r/2}[\epsilon^2]=1.$ 
\end{proof}

\begin{lem}
\label{lem.trunc_exp}
Let $\epsilon\sim \mathcal{N}(0,1).$ For any $t\geq 2$ and any $p\geq 0,$ $\E[|\epsilon|^p\mathbf{1}(|\epsilon|\geq t)]\leq (2+\sqrt{2/\pi}t^p) e^{-t^2/2}.$
\end{lem}

\begin{proof}
Notice that if this result is true for $p$ then also for all $p'\leq p.$ It is thus sufficient to prove the assertion for $p$ an even non-negative integer. Write  $F_p= \E[|\epsilon|^p\mathbf{1}(|\epsilon|\geq t)] = 2 \E[\epsilon^p\mathbf{1}(\epsilon\geq t)].$ If $p\geq 2,$ then, by partial integration, $F_p= \sqrt{2/\pi} t^{p-1} e^{-t^2/2} + F_{p-2}.$ By induction and using a simple bound for the geometric sum, we find for $t\geq 2,$ $F_p\leq \sqrt{2/\pi} t^p e^{-t^2/2} +F_0.$ The result follows from $F_0=2\P(\epsilon>t)\leq 2e^{-t^2}.$ 
\end{proof}

For the proof of the main theorem, an exponential inequality of the empirical $p$-th moment of i.i.d. Gaussian random variables for possibly large $p$ is required. For $p>2,$ exponential moments do not exist anymore and we can therefore not expect the usual $e^{-t^2}$ or $e^{-t}$ tail probability decay. Instead, the decay is $\exp(-t^{2/p}).$ An exponential inequality that comes arbitrarily close to the power $2/p$ can be found in Tao \cite{tao}, Exercise 2.1.7. Janson \cite{janson}, Theorem 6.7 and Theorem 6.12 shows that $\exp(-t^{2/p})$ is indeed the correct order.  For our approach, we also need explicit expressions of the constants, which cannot be deduced from the aforementioned results. Here, we give a new proof that allows to track the constants but is slightly suboptimal in the sense that we only obtain tail probability decay $\exp(-t^{2/(p+1)}).$

\begin{lem}
\label{lem.exp_ineq}
Let $\epsilon_k \sim \mathcal{N}(0,1),$ i.i.d. Then, for any $p\geq 1$ and any $t>0,$
\begin{align*}
	\P\big( \sum_{k=1}^d |\epsilon_k|^p - \E |\epsilon_k|^p \geq  6\sqrt{2} e p^{p/2} d^{1/2} t\big)\leq e^2\exp(-t^{2/(p+1)}).
\end{align*}
\end{lem}

\begin{proof}
For $t^{2/(p+1)}\leq 2$ there is nothing to show. We might therefore assume that $t^{2/(p+1)}>2.$ The Marcinkiewicz - Zygmund inequality by Ren and Liang \cite{ren} states that for independent, centered random variables $X_1, \ldots,X_d,$
\begin{align*}
	\E\big | \sum_{k=1}^d X_k \big| \leq (3\sqrt{2})^q q^{q/2} d^{q/2-1}\sum_{k=1}^d \E |X_k|^q, \quad q\geq 2.
\end{align*}
We have thus for the $q$-th moment $\E[| \sum_{k=1}^d |\epsilon_k|^p - \E |\epsilon_k|^p |^q] \leq (6\sqrt{2})^q q^{q/2} d^{q/2} \E[|\epsilon_1|^{pq}].$ With Lemma \ref{lem.mom_norm},  $\P( \sum_{k=1}^d |\epsilon_k|^p - \E |\epsilon_k|^p \geq 6\sqrt{2} e p^{p/2} d^{1/2} t)\leq q^{(p+1)q/2}e^{-q} t^{-q}.$ The result follows by choosing $q=t^{2/(p+1)}>2.$
\end{proof}

\begin{lem}
\label{lem.large_dvs_eps}
For any $j>J_0,$ any $p\geq 1,$ and $K_p$ as defined in \eqref{eq.Kp_def}, we have
\begin{align*}
	\sum_{k\in I_j} |\epsilon_{j,k}|^p < | I_j| K_p,
\end{align*}
with probability \ $\geq 1- e^2/n.$
\end{lem}

\begin{proof}
Use Lemma \ref{lem.exp_ineq} with $d=|I_j|$ and $t= \sqrt{|I_j|},$ Lemma \ref{lem.mom_norm} and $p+1\leq \log n/(\log \log n)^2,$ to obtain
\begin{align*}
		\P\big( \sum_{k=1}^d |\epsilon_k|^p \geq K_p |I_j| \big)\leq e^2e^{-|I_j|^{1/(p+1)}}
			\leq e^2 \exp(-|I_j|^{1/(p+1)})\leq e^2 \exp(-e^{(\log n) u_n /(p+1)})\leq \frac{e^2}n.
\end{align*}
\end{proof}

Define $\mT= \mT_n$ as the event
\begin{align}
	\max_{j\leq J,k \in I_j} |\epsilon_{j,k}| < 2\sqrt{\tfrac{\log n}n} \quad \text{and} \ \ 
	\sum_{k\in I_j} |\epsilon_{j,k}|^p < | I_j| K_p, \quad \text{for all} \ j\leq J, \ p\in \mP.
	\label{eq.mT_def}
\end{align}
Using a union bound, Lemma \ref{lem.large_dvs_eps} and the inequality $\P(|\epsilon|>t)\leq 2e^{-t^2/2}$ for $\epsilon\sim \mathcal{N}(0,1)$ give
\begin{align*}
	\P(\mT^c)\lesssim \frac{\log^2 n}{n}.
\end{align*}

\begin{lem}
\label{lem.norm_equiv}
On the event $\mT,$ for any $p\in \mP$ and for any $J_0<j\leq J,$
\begin{itemize}
\item[(i)] if $\max_{S: |S|=s}\sum_{k\in S} d_{j,k}^p\leq K_p n^{-p/2} |I_j|,$ then, $L_{j,p} >s,$
\item[(ii)] if for $S\subset I_j,$ $\sum_{k\in S} d_{j,k}^p\geq 4^{p} K_p n^{-p/2} |I_j|,$ then, $L_{j,p}\leq |S|,$
\item[(iii)] if $\max_{k\in I_j} |d_{j,k}| \leq \sqrt{\tfrac{\log n}n},$ then, $t_j= \max_{p\in \mP} t_{j,p},$
\item[(iv)] if $\max_{k\in I_j} |d_{j,k}| \geq 5\sqrt{\tfrac{\log n}n},$ then, $t_j=\infty.$
\end{itemize}
\end{lem}

\begin{proof}
On $\mT,$ the first statement follows due to
\begin{align*}
	\sum_{k\in S} Y_{j,k}^p \leq 2^{p-1} \sum_{k\in S} d_{j,k}^p + 2^{p-1} n^{-p/2} \sum_{k\in I_j} \epsilon_{j,k}^p
	< 2^p n^{-p/2} | I_j| K_p.
\end{align*}
The same argument can be used for the second assertion. $(iii)$ and $(iv)$ are a direct consequence from the first inequality in \eqref{eq.mT_def}.
\end{proof}

\begin{proof}[Proof of Proposition \ref{prop.properties_of_adapt_est}]
{\it (i):} Suppose that we are on the set $\mT.$ By Lemma \ref{lem.norm_equiv}, $\max_{k\in I_j} |d_{j,k}|\leq \sqrt{\tfrac{\log n}n}$ implies $t_j = \max_{p\in \mP} t_{j,p}.$ Let $p\in \mP$ be arbitrary. We show $L_{j,p}\geq 2.$ Suppose $L_{j,p}=1,$ then there exists $k^*,$ such that $|Y_{j,k^*}| \geq 2K_p^{1/p} n^{-1/2} |I_j|^{1/p}\geq 2\tfrac{\log n}{\sqrt{n}}\geq 3\sqrt{\tfrac{\log n}n}> \max_{k\in I_j} |Y_{j,k}| ,$ for $n>9.$  This is a contradiction and therefore, $L_{j,p}>1.$ Next, we show $L_{j,p}> K_pn^{-p/2} |I_j| Q^{-p}2^{jp(2\beta+1)}$ by deriving a contradiction. Suppose that $L_{j,p}\leq K_pn^{-p/2} |I_j| Q^{-p}2^{jp(2\beta+1)} \wedge |I_j|.$ For any subset $S\subset I_j$ with cardinality $L_{j,p},$
\begin{align*}
	\sum_{k\in S} |d_{j,k}|^p \leq L_{j,p} Q^p 2^{-\frac j 2 p(2\beta+1)}\leq K_p \frac{|I_j|}{n^{p/2}}
\end{align*}
and therefore Lemma \ref{lem.norm_equiv} (i) implies $L_{j,p}= |S| < L_{j,p},$ which is a contradiction. Therefore, we must have $L_{j,p}> K_pn^{-p/2} |I_j| Q^{-p}2^{\frac j 2  p(2\beta+1)}.$ Since also $L_{j,p}\geq 2,$ we find $L_{j,p}-1\geq L_{j,p}/2$ and $t_{j,p}\leq 4Q2^{-\frac j 2 (2\beta+1)}.$ Since $p\in \mP$ was arbitrary, also $t_j \leq 4Q2^{-\frac j 2 (2\beta+1)}.$ The assertion is proved since $\P(\mT^c)\lesssim \log^2 n/n.$

{\it (ii):} If $t_j=\infty$ then there is nothing to prove. Otherwise, it is enough to show $|\{k : |Y_{j,k}|>t_{j,p}\}|<L_{j,p}$ for all $p\in \mP.$ Suppose this is not true, that is, there exists $p^*\in \mP,$ such that $|\{k : |Y_{j,k}|>t_{j,p^*}\}|\geq L_{j,p^*}.$  Then, we can find a subset $R\subset \{k : |Y_{j,k}|>t_{j,p^*}\}$ of cardinality $L_{j,p^*}-1$ and
\begin{align*}
	\sum_{k\in R} Y_{j,k}^{p^*} > (L_{j,p^*}-1) t_{j,p^*}^{p^*} = 2^{p^*} K_{p^*} n^{-{p^*}/2} |I_j|.
\end{align*}
By the definition of $L_{j,p^*},$ we have $L_{j,p^*}-1= |R| \geq L_{j,p^*}.$ This is a contradiction and therefore, we must have $|\{k : |Y_{j,k}|>t_{j,p}\}|<L_{j,p}$ for all $p\in \mP.$
\end{proof}

\begin{rem}
\label{rem.prop_rem}
The proof implies that the set in Proposition \ref{prop.properties_of_adapt_est} (i) on which $t_j\leq 4Q2^{-\frac j2 (2\beta+1)}$ for all $J_0<j\leq J$ can be taken to be $\mT.$
\end{rem}

\begin{proof}[Proof of Theorem \ref{thm.non_adapt}]
Observe that for $j\leq J,$ $|\widetilde d_{j,k}-d_{j,k}| \leq |Y_{j,k}-d_{j,k}|\wedge 2Q2^{-\frac j2(2\beta+1)}.$

{\it $L^p$-risk for $1\leq p<\infty:$} Uniformly over $f\in \Theta(\beta,Q),$ using \eqref{eq.LP_embedding}, Jensen's inequality, and Lemma \ref{lem.mom_norm},
\begin{align*}
 \E\big[\|\widetilde f-f\|_p\big]
  &\lesssim \sum_{j=-1}^{J_n(\beta)}2^{j(\frac 12 -\frac 1p)} \E\big[\big(\sum_{k\in I_j}|\widetilde d_{j,k}-d_{j,k}|^p\big)^{1/p}\big]+ \sum_{j=J_n(\beta)+1}^\infty 2^{j(\frac 12 -\frac 1p)}  \big(\sum_{k\in I_j} |d_{j,k}|^p\big)^{1/p} \\
  &\leq \sum_{j=-1}^{J_n(\beta)} 2^{j(\frac 12 -\frac 1p)} \E^{1/p}\big[\sum_{k\in I_j} |\epsilon_{j,k}/\sqrt{n}|^p\big]+ O(n^{-\beta/(2\beta+1)}) \\
  &=O\big(\sqrt{\tfrac pn}2^{J_n(\beta)/2}+n^{-\beta/(2\beta+1)}\big)=O(\sqrt{p} n^{-\beta/(2\beta+1)}).
\end{align*}

{\it $L^\infty$-risk:} Let $\widetilde J_n(\beta)$ be an integer sequence such that $2^{\widetilde J_n(\beta)}\asymp (n/\log n)^{1/(2\beta+1)}.$ Due to \eqref{eq.LP_embedding}, the loss $\|\widetilde f-f\|_{\infty}$ can be bounded by a constant multiple of $  \sum_{j=-1}^{\infty} 2^{j/2} \max_k |\widetilde d_{j,k}-d_{j,k}|.$  Uniformly over $f\in \Theta(\beta,Q),$
\begin{align*}
  \sum_{j=-1}^{\infty} 2^{j/2} \max_k |\widetilde d_{j,k}-d_{j,k}|
  &\leq \sum_{j=-1}^{\widetilde J_n(\beta)} \frac{2^{j/2}}{\sqrt{n}} \max_k |\epsilon_{j,k}|+ \sum_{j=\widetilde J_n(\beta)+1}^{\infty} 2^{j/2} 2Q2^{-\frac j2(2\beta+1)}\\
  &\leq \sum_{j=-1}^{\widetilde J_n(\beta)} \sum_{k\in I_j} |\epsilon_{j,k}| \mathbf{1}( |\epsilon_{j,k}| \geq 2\sqrt{\log n})+O\big(\big(\tfrac{n}{\log n}\big)^{-\beta/(2\beta+1)}\big).
\end{align*}
Taking expectation, the result follows from Lemma \ref{lem.trunc_exp}.
\end{proof}

\begin{proof}[Proof of Theorem \ref{thm.main}]
Recall that $f\in \Theta(\beta,Q)$ implies $|d_{j,k}|\leq Q2^{-\frac j2(2\beta+1)}$ for all $j\geq -1,k\in I_j.$ Hence, there are integers $J_n(\beta)$ and $\widetilde J_n(\beta)$ such that  $J_0\leq \widetilde J_n(\beta)\leq J_n(\beta)\leq J,$
\begin{align}
  2^{\widetilde J_n(\beta)} \lesssim \big(\frac n{\log n}\big)^{1/(2\beta+1)}, \ \ 2^{J_n(\beta)} \asymp n^{1/2\beta+1}, \ \ \text{and} \  
  \sup_{f\in \Theta(\beta,Q)} \ \ \max_{j> \widetilde J_n(\beta),\ k\in I_j} |d_{j,k}|\leq  \sqrt{\tfrac{\log n}{n}}.
  \label{eq.tildeJnk_beta_conds}
\end{align}
All the estimates in this proof are uniformly over $f\in \Theta(\beta,Q)$ and $p\in \mP.$ In particular, $\lesssim$ means smaller or equal up to a constant multiple that only depends on $\beta,Q,$ and the underlying wavelet.

{\it $L^p$-risk for $1\leq p<\infty:$} Let $p_n := \log n/(\log \log n)^2-1.$ It is enough to consider $p\leq p_n/2$ since otherwise we can trivially bound the $L^p$-risk by the $L^\infty$-risk. For any $1\leq p\leq p_n/2,$ we can find a $p'\in \mP$ such that $p\leq p'< 2p.$ Since $\|\cdot \|_p\leq \|\cdot \|_{p'},$ it is thus enough to consider $p\in \mP.$ Using \eqref{eq.LP_embedding} and $\widehat d_{j,k}=0$ for $j>J,$ 
  \begin{align}
    \E\big[\|\widehat f-f\|_p\big]
    &\lesssim \E\Big[ \sum_{j=-1}^J 2^{j(\frac 12 -\frac 1p)}\big(\sum_{k\in I_j} |\widehat d_{j,k}-d_{j,k}|^p\big)^{1/p}\Big]  
    +   \sum_{j=J+1}^\infty 2^{j(\frac 12 -\frac 1p)}\big(\sum_{k\in I_j} |d_{j,k}|^p\big)^{1/p} \notag \\
    &=\E\Big[ \sum_{j=-1}^J 2^{j(\frac 12 -\frac 1p)}\big(\sum_{k\in I_j} |\widehat d_{j,k}-d_{j,k}|^p\big)^{1/p}\Big]  + O(n^{-\beta}).
    \label{eq.L2rate_first_decomp}
  \end{align}  
  
The main term $\sum_{j=-1}^J 2^{j(\frac 12 -\frac 1p)}\big(\sum_{k\in I_j} |\widehat d_{j,k}-d_{j,k}|^p\big)^{1/p}$ can be bounded from above by the sum of the five terms 
%to split (III) and (IV) we used that (a+b)^{1/p}\leq a^{1/p}+b^{1/p}
\begin{align*}
  (I)&:= \sum_{j=J_0+1}^J 2^{j(\frac 12 -\frac 1p)}\big(\sum_{k\in I_j} |\widehat d_{j,k}-d_{j,k}|^p\big)^{1/p}\mathbf{1}(\mathcal{T}^c) \\
  (II)&:= \sum_{j=J_0+1}^{J_n(\beta)} 2^{j(\frac 12 -\frac 1p)}\big(\sum_{k\in I_j} |d_{j,k}|^p\big)^{1/p} \mathbf{1}(\mathcal{T} \cap\{t_j=0\}) \\
  (III)&:= \sum_{j=-1}^{J_n(\beta)} 2^{j(\frac 12 -\frac 1p)}\big(\sum_{k\in I_j} |Y_{j,k}- d_{j,k}|^p\big)^{1/p}   \\
  (IV)&:= 2\sum_{j=J_0+1}^{J_n(\beta)} 2^{j(\frac 12 -\frac 1p)}\big(\sum_{k\in I_j: |Y_{j,k}|>t_j} t_j^p + |d_{j,k}|^p\big)^{1/p} \mathbf{1}(\mathcal{T} \cap\{t_j>0\} )\\
  (V)&:= \sum_{j=J_n(\beta)+1}^J2^{j(\frac 12 -\frac 1p)}\big(\sum_{k\in I_j} |\widehat d_{j,k}-d_{j,k}|^p\big)^{1/p}\mathbf{1}(\mathcal{T}).
\end{align*}
In the following we bound the expectation of $(I)-(V),$ separately.

{\it (I):} By definition $|\widehat{d}_{j,k}|\leq |Y_{j,k}|$ and therefore 
\begin{align*}
	(\widehat d_{j,k}-d_{j,k})^p
	&\leq 2^{p-1}\widehat d_{j,k}^p+2^{p-1}d_{j,k}^p \\
	&\leq 4^{p}d_{j,k}^p+2^p n^{-p/2} \epsilon_{j,k}^p \\
	&\leq 4^{p}d_{j,k}^p+2^p n^{-p/2} (2p\log n)^{p/2}+ 2^p n^{-p/2} \epsilon_{j,k}^p\mathbf{1}(|\epsilon_{j,k}|>\sqrt{2p\log n}).
\end{align*}
Thus, we find
\begin{align*}
	(I)\lesssim \sqrt{p\log n}  \cdot \mathbf{1}(\mathcal{T}^c) + \frac{2}{\sqrt{n}} \sum_{j=J_0+1}^J 2^{j(\frac 12 -\frac 1p)}\big(\sum_{k\in I_j} |\epsilon_{j,k}|^p\mathbf{1}(|\epsilon_{j,k}|>\sqrt{2p\log n})\big)^{1/p}
\end{align*}
Taking expectation, applying $\P(\mT^c)\lesssim \log^2 n/n,$  $E[|X|^{1/p}]\leq E^{1/p}[|X|],$ and Lemma \ref{lem.trunc_exp} with $t=\sqrt{2p\log n},$
\begin{align*}
	\E\big[\sum_{j=J_0+1}^J 2^{j(\frac 12 -\frac 1p)}\big(\sum_{k\in I_j} |\widehat d_{j,k}-d_{j,k}|^p\big)^{1/p}\mathbf{1}(\mathcal{T}^c)\big]
	\lesssim \sqrt{p} \frac{\log^{5/2} n}{n} + \frac{\sqrt{p\log n}}{n}.
\end{align*}

{\it (II):} If $t_j=0,$ then also $L_{j,p}=\infty$ for all $p\in \mP.$  By Lemma \ref{lem.norm_equiv} (ii), we must consequently have that $\sum_{k\in I_j} d_{j,k}^p \leq 4^p K_p n^{-p/2} |I_j|,$ on $\mT.$ Thus,
\begin{align*}
	(II)\leq \frac{4}{\sqrt{n}} K_p^{1/p} \sum_{j=J_0+1}^{J_n(\beta)} 2^{j(\frac 12 -\frac 1p)} |I_j|^{1/p}
	\lesssim \sqrt{p} n^{-1/2} 2^{J_n(\beta)/2}\lesssim \sqrt{p} n^{-\beta/(2\beta+1)}.
\end{align*}

{\it (III):} Lemma \ref{lem.mom_norm} gives
\begin{align*}
	\E\big[\sum_{j=-1}^{J_n(\beta)} 2^{j(\frac 12 -\frac 1p)}\big(\sum_{k\in I_j} |Y_{j,k}- d_{j,k}|^p\big)^{1/p} \big]
	\lesssim \sum_{j=-1}^{J_n(\beta)} 2^{j(\frac 12 -\frac 1p)} |I_j|^{1/p} n^{-1/2}\E^{1/p}[\epsilon^{p}] \lesssim \sqrt{p} n^{-\beta/(2\beta+1)}.
\end{align*}

{\it (IV):} For any $J_0< j\leq J,$ let $R_j = \{k\in I_j : |Y_{j,k}|>t_j\}.$ By Proposition \ref{prop.properties_of_adapt_est} (ii), $R_j<\min_{p\in \mP} L_{j,p}.$ Suppose that there exists a $p\in \mP$ with $\sum_{k\in R_j} d_{j,k}^p \geq 4^p K_p n^{-p/2} |I_j|.$ Then, on $\mT,$ Lemma \ref{lem.norm_equiv} (ii) implies $|R_j|\geq L_{j,p}$ which is a contradiction since $R_j<\min_{p\in \mP} L_{j,p}.$ Thus, for any $p\in \mP,$ $\sum_{k\in R_j} d_{j,k}^p \leq 4^p K_p n^{-p/2} |I_j|.$ This allows us to bound
\begin{align*}
	(IV) \leq 12 \sum_{j=J_0+1}^{J_n(\beta)} 2^{j(\frac 12 -\frac 1p)} K_p^{1/p} n^{-1/2} |I_j|^{1/p} \lesssim \sqrt{p} n^{-\beta/(2\beta+1)}.
\end{align*}

{\it (V):} Since $\widetilde J_n(\beta)\leq J_n(\beta),$ we have by \eqref{eq.tildeJnk_beta_conds}, $\max_{k\in I_j} |d_{j,k}|\leq \sqrt{\log n/n}$ for all $J_n(\beta)\leq j\leq J.$ Using Proposition \ref{prop.properties_of_adapt_est} (i) and Remark \ref{rem.prop_rem}, $t_j\leq 4Q2^{-\frac j2 (2\beta+1)}$ on $\mT$ and therefore, $(V)\leq 5Q \sum_{j=J_n(\beta)}^J  2^{j(\frac 12 -\frac 1p)} 2^{-\frac j2 (2\beta+1)} |I_j|^{1/p}\lesssim n^{1/(2\beta+1)}.$

The estimates $(I)-(V)$ together with \eqref{eq.L2rate_first_decomp} show that for any $p\in \mP$ the $L^p$-risk is bounded by const.$ \times \sqrt{p} n^{-\beta/(2\beta+1)}$ and this completes the proof for the adaptive rate under $L^p$-loss.

{\it $L^\infty$-risk:} Thanks to the embedding \eqref{eq.LP_embedding}, we can bound the ${L^\infty}$-norm of the estimator $\widehat f$ by a constant multiple of $ \sum_{j=-1}^{J} 2^{j/2} \max_k |\widehat d_{j,k}-d_{j,k}|+\sum_{j=J+1}^\infty 2^{j/2}\max_k |d_{j,k}|.$ The second term is of the (negligible) order $n^{-\beta}.$ It thus remains to show
\begin{align}
  \E\big[\sum_{j=-1}^{J} 2^{j/2} \max_k |\widehat d_{j,k}-d_{j,k}|\big]\lesssim \big(\tfrac n{\log n}\big)^{-\beta/(2\beta+1)}.
  \label{eq.Linf_to_show_mainTHM}
\end{align}
In the following we bound the expectation of the three summands
\begin{align*}
  (i) &:= \sum_{j=-1}^{J} 2^{j/2} \max_k |\widehat d_{j,k}-d_{j,k}|\mathbf{1}(\mathcal{T}^c) ,\notag \\
  (ii) &:= \sum_{j=-1}^{\widetilde J_n(\beta)} 2^{j/2} \max_k |\widehat d_{j,k}-d_{j,k}| \mathbf{1}(\mathcal{T}), \\
  (iii) &:=\sum_{j=\widetilde J_n(\beta)+1}^{J} 2^{j/2} \max_k |\widehat d_{j,k}-d_{j,k}| \mathbf{1}(\mathcal{T}). 
\end{align*}

{\it (i):} Using $|\widehat d_{j,k}|\leq |d_{j,k}|+ n^{-1/2} |\epsilon_{j,k}|\leq |d_{j,k}|+ 2\sqrt{\log n/n} + n^{-1/2} |\epsilon_{j,k}|\mathbf{1}(|\epsilon_{j,k}|\geq 2\sqrt{\log n}) $ and Lemma \ref{lem.trunc_exp},
\begin{align*}
  \E \big[\sum_{j=-1}^{J} 2^{j/2} \max_{k \in I_j} |\widehat d_{j,k}-d_{j,k}|\mathbf{1}(\mathcal{T}^c)\big]
  &\lesssim (\sqrt{\log n}) \P(\mathcal{T}^c)+\E\big[\sum_{j=-1}^{J} \sum_{k \in I_j} |\epsilon_{j,k}|\mathbf{1}(|\epsilon_{j,k}|\geq 2\sqrt{\log n})\big]  \\  
  &= O\big(\tfrac{\log^{5/2} n}{n}\big).
\end{align*}

{\it (ii):} Work on $\mT$ and let $-1\leq j\leq  \widetilde J_n(\beta)$ be arbitrary. We need to consider $I_j' =\{k \in I_j : \widehat{d}_{j,k} = Y_{j,k}\}$ and $I_j'' =\{k \in I_j : \widehat{d}_{j,k} = t_j\},$ separately. Notice that $  2^{j/2}\max_{k\in I_j'}|\widehat d_{j,k}-d_{j,k}|\mathbf{1}(\mathcal{T})
  \leq 2^{j/2}n^{-1/2} \max_k  |\epsilon_{j,k}| \mathbf{1}(\mathcal{T})
  \leq 2^{j/2} 2\sqrt{\log n/n}.$ If $k\in I_j'',$ then $t_j<\infty$ and  $t_j = \max_{p\in \mP} t_{j,p}.$ In a first step, we  show that $t_j\leq 20\sqrt{\log n/n}.$ To see this recall that by Lemma \ref{lem.norm_equiv} (iv),  $\max_k |d_{j,k}|\leq 5\sqrt{\log n/n}.$ Fix $p\in \mP.$ For any subset $R\subset I_k$ with $|R|\leq K_p 5^{-p}|I_j| /\log ^{p/2} n,$ we must have $\sum_{k\in R} d_{j,k}^p \leq K_p n^{-p/2} |I_j|.$ Applying Lemma \ref{lem.norm_equiv} (i) shows that $L_{j,p} > K_p 5^{-p} |I_j|/\log ^{p/2} n .$ The r.h.s. is larger than one, because of $j\geq J_0,$ and $p+1 \leq \log n/(\log \log n)^2.$ Consequently, $L_{j,p}-1 > \tfrac 12  K_p 5^{-p} |I_j|/\log ^{p/2} n $ which implies that $t_{j,p}\leq 20\sqrt{\log n/n}.$ Since this holds for any $p\in \mP,$ also $t_j\leq 20\sqrt{\log n/n}$ and
\begin{align*}
 2^{j/2}\max_{k\in I_j''} |\widehat d_{j,k}-d_{j,k}| \mathbf{1}(\mathcal{T})
  \leq 2^{j/2} (t_j+ |d_{j,k}|)\mathbf{1}(\mathcal{T})
  \leq 25 \cdot 2^{j/2} \sqrt{\frac{\log n}n}.
\end{align*}
Recall that $j$ in $J_0\leq j\leq \widetilde J_n(\beta)$ was arbitrary. Combining the upper bounds for $k\in I_j'$ and $k\in I_j''$ gives
\begin{align*}
  (ii) \lesssim \sum_{j=-1}^{\widetilde J_n(\beta)} 2^{j/2} \sqrt{\frac{\log n}n}
  \lesssim \big(\frac n{\log n}\big)^{-\beta/(2\beta+1)}.
\end{align*}

{\it (iii):} Since $\max_{k\in I_j}|d_{j,k}|\leq  \sqrt{\log n/n},$ we have by Proposition \ref{prop.properties_of_adapt_est} (i) and Remark \ref{rem.prop_rem} that $|\widehat{d_{j,k}}|\leq t_j \leq 4Q2^{-\frac j2(2\beta+1)},$ on the event $\mT.$ This proves $(iii) \leq  5Q\sum_{j=\widetilde J_n(\beta)+1}^{J} 2^{- j\beta} \lesssim (n/\log n)^{-\beta/(2\beta+1)}.$

Combining the upper bounds for $(i)-(iii)$ shows \eqref{eq.Linf_to_show_mainTHM} and this completes the proof for the $L^\infty$-risk.
\end{proof}

\begin{proof}[Proof of Theorem \ref{thm.neccesary_cond}]
Let $\beta, Q\in (0,\infty)$ be arbitrary. Notice that $f_j:= \sum_{k\in I_j} Q2^{-\frac j2(2\beta+1)} \psi_{j,k} \newline \in \Theta(\beta,Q).$ Choose an integer $J_n=J_n(\beta),$ such that $n^{1/(2\beta+1)} \lesssim 2^{J_n}\leq (Q^2n/(16K_2))^{1/(2\beta+1)}$ and consider $f^0:= f_{J_n}\in \Theta(\beta,Q).$ By Lemma \ref{lem.norm_equiv}, we find $t_j\geq t_{j,2}>0$ on $\mT$ and thus,
\begin{align}
  \sup_{f\in \Theta(\beta,Q)} \E\big[\|\widehat f_1-f\|_{\infty}\big]
  \geq \E\big[\|\widehat f_1-f^0\|_{\infty} \mathbf{1}(\mT)\big]
  = \frac 1{\sqrt{n}} \E\big[\|\sum_{k\in I_{J_n}} \epsilon_{J_n,k} \psi_{J_n,k}\|_{\infty} \mathbf{1}(\mT)\big].
  \label{eq.nec_thm_eval_at_funct}
\end{align}
Thanks to the compact support of the wavelet function, we can always pick a subset $I\subset I_{J_n}$ of cardinality $\gtrsim 2^{J_n}$ such that for any $k_1, k_2 \in I,$ $k_1\neq k_2,$ the wavelet functions $\psi_{J_n,k_1} $ and $\psi_{J_n,k_2}$ have disjoint support. Define $k^*$ by
\begin{align*}
 k^*\in \argmax_{k\in I}\ \epsilon_{J_n,k} 
\end{align*}
and consider the sets $\mathcal{R}:= \{\epsilon_{J_n,k^*}\geq \sqrt{\log |I|}\}$ and $U(k^*):=\{k: k\neq k^*$ and the support of $\psi_{J_n,k}$ and the support of $\psi_{J_n,k^*}$ have non-empty intersection $\}.$ It follows from the extreme value behavior of standard normal random variables (cf. Embrechts et {\it al.} \cite{emb}, p.145) that $\P(\mathcal{R}^c)< 1/2$ for sufficiently large $n.$ For sufficiently small $\delta>0,$ there exists a random sequence $x_n\in (0,1)$ such that $$\inf_n |2^{-J_n/2}\psi_{J_n,k^*}(x_n)|\geq \delta>0.$$ By triangle inequality,
\begin{align}
  &\|\sum_{k \in I_{J_n}} \epsilon_{J_n,k} \psi_{J_n,k}\|_{\infty} \mathbf{1}(\mT) \geq |\epsilon_{J_n,k^*} \psi_{J_n,k^*}(x_n)| \mathbf{1}(\mT \cap \mathcal{R})
  - |\sum_{k\in U(k^*)} \epsilon_{J_n,k} \psi_{J_n,k}(x_n)|.
  \label{eq.nec_thm_lb_main_term}
\end{align}
To bound the second term, notice that $U(k^*)$ and $I$ are disjoint by construction. Since $k^*$ is a function of $\{\epsilon_{j,k}: k\in I\}$ it is independent of $\epsilon_{J_n,k}$ for all $k\in U(k^*).$ The same holds for $x_n.$ Using $\E[\cdot]=\E[\E[\cdot|k^*,x_n]],$ it follows
\begin{align}
 \E\big[|\sum_{k\in U(k^*)} \epsilon_{J_n,k} \psi_{J_n,k}(x_n)|\big]
  \leq 2^{J_n/2} \|\psi\|_{\infty} \E\big[\sum_{k\in U(k^*)} |\epsilon_{J_n,k}|\big]
  \lesssim 2^{J_n/2}.
  \label{eq.nec_thm_remainder_ub}
\end{align}
In order to find a lower bound of the first term in \eqref{eq.nec_thm_lb_main_term}, observe that $\P(\mathcal{T}\cap \mathcal{R}) \geq 1-\P(\mathcal{T}^c)-\P(\mathcal{R}^c)\geq 1/4$ for sufficiently large $n,$ since $\P(\mathcal{T}^c)\rightarrow 0$ and $\P(\mathcal{R}^c)<1/2.$ Therefore,
\begin{align*}
  \E\big[|\epsilon_{J_n,k^*} \psi_{J_n,k^*}(x_n)| \mathbf{1}(\mT \cap \mathcal{R})\big]
  \geq \frac {\delta}4 2^{J_n/2} \sqrt{\log |I|}.  
\end{align*}
Together with \eqref{eq.nec_thm_eval_at_funct}, \eqref{eq.nec_thm_lb_main_term}, and \eqref{eq.nec_thm_remainder_ub}  the result follows since by definition of $I,J_n,$
\begin{align*}
  n^{-1/2} 2^{J_n/2} \sqrt{\log |I|}
  \gtrsim \big(\frac n{\log n} \big)^{-\beta/(2\beta+1)} (\log n)^{1/(4\beta+2)}.
\end{align*}
\end{proof}

\bibliographystyle{acm}    %acm   % (uses file "plain.bst")
\bibliography{refsPart1}           % expects file "refsPart1.bib"

\end{document}